\documentclass{birkart}
\usepackage{hyperref}

\newcommand{\arxiv}[1]{\href{http://arxiv.org/abs/#1}{arXiv:#1}}
\newcommand*{\mailto}[1]{\href{mailto:#1}{\nolinkurl{#1}}} 

\newtheorem{theorem}{Theorem}[section]
\newtheorem{lemma}[theorem]{Lemma}
\newtheorem{corollary}[theorem]{Corollary}
\newtheorem{remark}[theorem]{Remark}

\newcommand{\R}{{\mathbb R}}
\newcommand{\N}{{\mathbb N}}
\newcommand{\Z}{{\mathbb Z}}
\newcommand{\C}{{\mathbb C}}
\newcommand{\M}{{\mathbb M}}

\newcommand{\nn}{\nonumber}
\newcommand{\be}{\begin{equation}}
\newcommand{\ee}{\end{equation}}
\newcommand{\bea}{\begin{eqnarray}}
\newcommand{\eea}{\end{eqnarray}}

\newcommand{\ol}{\overline}
\newcommand{\ti}{\tilde}
\newcommand{\what}{\widehat}

\newcommand{\id}{{\rm 1\hspace{-0.6ex}l}}
\newcommand{\E}{\mathrm{e}}
\newcommand{\I}{\mathrm{i}}

\newcommand{\KdV}{\mathrm{KdV}}
\newcommand{\tr}{\mathrm{tr}}
\newcommand{\im}{\mathrm{Im}}
\newcommand{\re}{\mathrm{Re}}

\newcommand{\di}{\mathcal{D}}

\newcommand{\dirho}{\di_{\rho}}
\newcommand{\dirhos}{\di_{\rho}^*}

\newcommand{\Rg}[1]{R_{2g+1}^{1/2}(#1)}
\newcommand{\Rgo}{R_{2g+1}^{1/2}}

\newcommand{\sig}{\sigma}
\newcommand{\lam}{\lambda}

\newcommand{\om}{\omega}

\renewcommand{\labelenumi}{(\roman{enumi})}
\numberwithin{equation}{section}

\pagespan{00}{99}

\begin{document}

\title[Trace Formulas for Schr\"odinger Operators]{Trace Formulas for
Schr\"odinger Operators in Connection with Scattering Theory for Finite-Gap Backgrounds}

\author[A.\ Mikikits-Leitner]{Alice Mikikits-Leitner}
\address{Faculty of Mathematics\\ University of Vienna\\
Nordbergstrasse 15\\ 1090 Wien\\ Austria}
\email{\mailto{Alice.Mikikits-Leitner@univie.ac.at}}
\urladdr{\url{http://www.mat.univie.ac.at/~alice/}}

\author[G.\ Teschl]{Gerald Teschl}
\address{Faculty of Mathematics\\ University of Vienna\\
Nordbergstrasse 15\\ 1090 Wien\\ Austria\\ and International Erwin Schr\"odinger
Institute for Mathematical Physics\\ Boltzmanngasse 9\\ 1090 Wien\\ Austria}
\email{\mailto{Gerald.Teschl@univie.ac.at}}
\urladdr{\url{http://www.mat.univie.ac.at/~gerald/}}

\thanks{Work supported by the Austrian Science Fund (FWF) under Grant No.\ Y330}
\thanks{in Spectral Theory and Analysis, J. Janas (ed.) et al., 107--124, Oper. Theory Adv. Appl. {\bf 214}, Birkh\"auser, Basel, 2011}


\begin{abstract}
We investigate trace formulas for one-dimensional Schr\"odinger operators which are
trace class perturbations of quasi-periodic finite-gap operators using Krein's
spectral shift theory. In particular, we establish the conserved quantities for the
solutions of the Korteweg--de Vries hierarchy in this class and relate them to
the reflection coefficients via Abelian integrals on the underlying hyperelliptic
Riemann surface.
\end{abstract}

\maketitle

\section{Introduction}

Trace formulas for one-dimensional (discrete and continuous) Schr\"odinger operators
have attracted an enormous amount of interest recently (see e.g.\ \cite{dk}, \cite{glmz}, \cite{ks}, \cite{lns},
\cite{npy}, \cite{sizl}, \cite{zl}). However, most  results are in connection with scattering theory
for a constant background. On the other hand, scattering theory for one-dimensional
Schr\"odinger operators with periodic background is a much older topic first investigated by Firsova
in a series of papers \cite{F1}--\cite{F3}. Nevertheless, many questions which have long been
answered in the constant background case are still open in this more general setting.

The aim of the present paper is to help filling some of these gaps. To this end, we want to find the analog
of the classical trace formulas in scattering theory for the case of a quasi-periodic, finite-gap background. In the case of
zero background it is well-known that the transmission coefficient is the perturbation determinant in the
sense of Krein \cite{krein} (see e.g., \cite{jp}, \cite{sres}, \cite{yafams} see also \cite{GM}, \cite{gps} and the references
therein for generalizations to non trace class situations) and our first aim is to establish this result for the case
considered here; thereby establishing the connection with
Krein's spectral shift theory. Our second aim is to find a representation of the transmission
coefficient in terms of the scattering data --- the analog of the classical Poisson--Jensen formula.

Moreover, scattering theory for one-dimensional Schr\"odinger operators is not only interesting
in its own right, it also constitutes the main ingredient of the inverse scattering transform for
the Korteweg--de Vries (KdV) hierarchy (see, e.g., \cite{evh}, \cite{M}). Again the case of decaying solutions is classical
and trace formulas for this case were studied exhaustively in the past (cf.\ \cite{conl} and the references therein).
Here we want to investigate the case of Schwartz type perturbations of a given finite-gap solution.
The Cauchy problem for the KdV equation with initial conditions in this class was only solved recently by
Egorova, Grunert, and Teschl \cite{EGT} (see also \cite{ET}, \cite{ET2}, \cite{F4}).
Since the transmission coefficient is invariant when our Schr\"odinger
operator evolves in time with respect to some equation of the KdV hierarchy, the
corresponding trace formulas provide the conserved quantities for the
KdV hierarchy in this setting.

Our work extends previous results for Jacobi operators by Michor and Teschl \cite{mtqptr}, \cite{tag}.
For trace formulas in the pure finite-gap case see Gesztesy, Ratnaseelan, and Teschl \cite{GRT} and
Gesztesy and Holden \cite{GH}.

\section{Notation}

We assume that the reader is familiar with quasi-periodic, finite-gap Schr\"odinger operators which
arise naturally as the stationary solutions of the KdV hierarchy.
Hence we only briefly recall some notation and refer to the monograph
\cite{GH} (see also \cite{M}) for further information.

Let
\be
H_q = -\frac{d^2}{dx^2}+V_q(x)
\ee
be a finite-gap Schr\"odinger operator in $L^2(\R )$ whose spectrum consists of $g+1$ bands:
\be
\sig(H_q) = \bigcup_{j=0}^{g-1} [E_{2j},E_{2j+1}]\cup [E_{2g},\infty).
\ee
It is well-known that $H_q$ is associated with the Riemann surface $\M$
of the function
\be
\Rg{z}, \qquad R_{2g+1}(z) = \prod_{j=0}^{2g} (z-E_j), \qquad
E_0 < E_1 < \cdots < E_{2g},
\ee
$g\in \N_0$. $\M$ is a compact, hyperelliptic Riemann surface of genus $g$.
Here $\Rg{z}$ is chosen to have branch cuts along the spectrum with the sign fixed
by the asymptotic behavior $\Rg{z}=\sqrt{z}z^g+\dots$ as $z\to \infty$.

A point on $\M$ is denoted by $p = (z, \pm \Rg{z}) = (z, \pm)$, $z \in \C$.
The point at infinity is denoted by $p _\infty =(\infty, \infty)$.
We use $\pi(p) = z$ for the projection onto the extended complex plane
$\C \cup \{\infty\}$.  The points $\{(E_{j}, 0), 0 \leq j \leq 2 g\}  \cup \{p_\infty\}\subseteq \M$ are 
called branch points and the sets 
\be
\Pi_\pm = \{ (z, \pm \Rg{z}) \mid z \in \C\backslash \Sigma\} \subset \M,
\qquad \Sigma= \sig(H_q),
\ee
are called upper and lower sheet, respectively. Note that the boundary of
$\Pi_\pm$ consists of two copies of $\Sigma$ corresponding to the
two limits from the upper and lower half plane.

We recall that upon fixing the spectrum $\sig(H_q)$, the operator $H_q$ is uniquely
defined by choosing a Dirichlet divisor
\be\label{divisor}
\left\{\big(\mu_1,\sig_1\big), \dots, \big(\mu_g,\sig_g\big)\right\}, \qquad \mu_j \in [E_{2j-1},E_{2j}].
\ee
For every $z\in\C$ the Baker--Akhiezer functions $\psi_{q,\pm}(z,x)$ are
two (weak) solutions of $H_q \psi = z \psi$. They are the two branches of one function
which is meromorphic on $\M\backslash\{p_\infty\}$ with simple poles at the Dirichlet
divisor \eqref{divisor} and simple zeros at some other points
\be
\left\{\big(\mu_1(x),\sig_1(x)\big), \dots, \big(\mu_g(x),\sig_g(x)\big)\right\}, \qquad \mu_j(x) \in [E_{2j-1},E_{2j}],
\ee
which can be computed from the Dubrovin equations
\be
\mu_j'(x) = \frac{-2 \sig_j(x)\Rg{\mu_j(x)}}{\prod_{k\ne j} \mu_j(x)-\mu_k(x)}
\ee
using the initial conditions $\mu_j(0) =\mu_j$, $1\le j \le g$. Moreover, $V_q(x)$ is explicitly given by
the trace formula
\be\label{trVq}
V_q(x)=E_0+\sum_{j=1}^g\big(E_{2j-1}+E_{2j}-2\mu_j(x)\big).
\ee
The Baker--Akhiezer functions are linearly independent away
from the band-edges $\{E_j\}_{j=0}^{2g}$ since their Wronskian is given by
\be
W(\psi_{q,-}(z), \psi_{q,+}(z)) = \frac{2\I \Rg{z}}{\prod_{j=1}^g (z-\mu_j)}.
\ee
Here $W_x(f,g)=f(x)g'(x)-f'(x)g(x)$ denotes the usual Wronskian and $\mu_j$ are
the Dirichlet eigenvalues at base point $x_0=0$.
We recall that $\psi_{q,\pm}(z,x)$ have the form
\be\label{decomppsipm}
\psi_{q,\pm}(z,x) =\theta_{q,\pm}(z,x) \exp (\pm \I xk(z)),
\ee
where $\theta_{q,\pm}(z,x)$ is quasi-periodic with respect to $x$ and
\be
k(z)=-\int_{E_0}^{p}\omega_{p_{\infty},0}, \qquad p=(z,+),
\ee
denotes the quasi\-momentum map. Here $\omega_{p_{\infty},k}$ is a normalized
Abelian differential of the second kind with a single pole at $p_{\infty}=(\infty,\infty)$ and principal part
$\zeta^{-k-2}d\zeta$ where $\zeta=z^{-1/2}$. It is explicitly given by
\be
\omega_{p_{\infty},0}=-\frac{\prod_{j=1}^g(\pi-\lambda_j)d\pi}{2 \Rgo},
\ee
where $\lambda_j\in (E_{2j-1},E_{2j})$, $1\leq j\leq g$.
In particular, $\big|\E ^{\I k(z)}\big|<1$ for
$z\in\C\backslash\sig(H_q)$ and $|\E ^{\I k(z)}|=1$ for $z\in\sig(H_q)$.

\section{Asymptotics of Jost solutions}
\label{secJS}

After we have these preparations out of our way, we come to
the study of short-range perturbations $H$ of $H_q$ associated with a potential $V$
satisfying $V(x) \rightarrow V_q(x)$ as $|x|\rightarrow \infty$.  More precisely, we will
make the following assumption throughout this paper:

Let
\be
H = -\frac{d^2}{dx^2}+V(x)
\ee
be a perturbation of $H_q$ such that
\be                         \label{hypo}
\int_{-\infty}^{+\infty} \big|V(x) - V_q(x)\big| dx <\infty.
\ee
We first establish existence of Jost solutions, that is, solutions of the
perturbed operator which asymptotically look like the Baker--Akhiezer solutions.

\begin{theorem}
Assume \eqref{hypo}. For every $z \in \C\backslash\{E_j\}_{j=0}^{2g}$ there exist (weak) solutions $\psi_\pm(z, .)$
of $H \psi = z \psi$ satisfying
\be \label{jost1}
\lim_{x \rightarrow \pm \infty}
\E^{\mp \I xk(z)} \big( \psi_\pm(z, x) - \psi_{q,\pm}(z, x) \big) = 0,
\ee
where $\psi_{q,\pm}(z, .)$ are the Baker--Akhiezer functions.
Moreover, $\psi_\pm(z, .)$ are continuous 
(resp.\ holomorphic) with respect to $z$ whenever $\psi_{q,\pm}(z, .)$ are and
\be \label{jost2}
\big| \E^{\mp \I xk(z)} \big( \psi_\pm(z, x) - \psi_{q,\pm}(z, x) \big)\big|\leq C(z),
\ee
where $C(z)$ denotes some constant depending only on $z$.
\end{theorem}
\begin{proof}
Since $H\psi=z\psi$ is equivalent to $(H_q-z)\psi=-\what{V}\psi$, where $\what{V}=V-V_q$,
we can use the variation of constants formula to obtain the usual Volterra integral equations for the Jost functions,
\begin{align} \label{asyjost}
\psi_\pm(z,x)=\psi_{q,\pm}(z,x)-&\frac{1}{W(\psi_{q,+},\psi_{q,-})}\int_{x}^{\pm \infty} \big( \psi_{q,-}(z,x)\psi_{q,+}(z,y)-\nn\\
-&\psi_{q,-}(z,y)\psi_{q,+}(z,x)\big)\what{V}(y)\psi_\pm(z,y)dy.
\end{align}
Moreover, introducing $\ti\psi_\pm(z,x) = \E^{\mp \I xk(z)} \psi_\pm(z, x)$ the resulting integral equation can be solved using
the method of successive iterations in the usual way. This proves the claims.
\end{proof}

\begin{theorem} \label{thmjost}
Assume \eqref{hypo}. The Jost functions have the following asymptotic behavior
\be                 \label{B4jost}
\psi_\pm(z,x) =  \psi_{q,\pm}(z,x)\Big( 1\mp\frac{1}{2\I \sqrt{z}} \int_{x}^{\pm \infty}\big(V(y)-V_q(y)\big) dy+o(z^{-1/2})\Big),
\ee
as $z\to \infty$, with the error being uniformly in $x$.
\end{theorem}

\begin{proof}
Invoking \eqref{asyjost} we have
\begin{align} \label{Volterrafrac}
\frac{\psi_\pm (z,x)}{\psi_{q,\pm}(z,x)}
=1&-\frac{1}{W( \psi_{q,+},\psi_{q,-})}\int_{x}^{\pm \infty}\Big( \psi _{q,-}(z,x)\psi _{q,+}(z,y)\frac{\psi_{q,\pm}(z,y)}{\psi_{q,\pm}(z,x)}\nn \\
&- \psi _{q,-}(z,y)\psi _{q,+}(z,x)\frac{\psi_{q,\pm}(z,y)}{\psi_{q,\pm}(z,x)}\Big) \what{V}(y)\frac{\psi_\pm(z,y)}{\psi_{q,\pm}(z,y)}dy\nn \\
=1& \mp \int_{x}^{\pm \infty}\Big( G_q(z,x,x)\frac{\psi_{q,\pm}(z,y)^2}{\psi_{q,\pm}(z,x)^2}-G_q(z,y,y)\Big)\what{V}(y)\frac{\psi_\pm(z,y)}{\psi_{q,\pm}(z,y)}dy,
\end{align}
where 
\be \label{defGq}
G_q(z,x,y)=\frac{1}{W( \psi_{q,+},\psi_{q,-})}\left\{ \begin{array}{cc}
\psi_{q,+}(z,x)\psi_{q,-}(z,y), & x\geq y, \\ \psi_{q,+}(z,y)\psi_{q,-}(z,x), & x\leq y,
\end{array}
\right.
\ee
is the Green function of $H_q$.
We have
\be
G_q(z,x,x)=\frac{\psi_{q,+}(z,x)\psi_{q,-}(z,x)}{W(\psi_{q,+},\psi_{q,-})} =\frac{\I \prod_{j=1}^g\big(z-\mu_j(x)\big)}{2\Rg{z}}.
\ee
Hence for $z$ near $\infty$ one infers
\be \label{asyGreen}
G_q(z,x,x)=\frac{\I }{2\sqrt{z}} \Big( 1+\frac{1}{2}V_q(x)\frac{1}{z}+O\big( \frac{1}{z^2}\big) \Big),
\ee
where we made use of the trace formula \eqref{trVq}.
Next we insert \eqref{asyGreen} into \eqref{Volterrafrac} such that iteration implies
\be
\frac{\psi_\pm(z,x)}{\psi_{q,\pm}(z,x)}
=1\mp\frac{\I }{2\sqrt{z}}\Big( \int_x^{\pm \infty}\frac{\psi_{q,\pm}(z,y)^2}{\psi_{q,\pm}(z,x)^2}\what{V}(y)dy-
\int_{x}^{\pm \infty}\what{V}(y)dy\Big)+O\big(\frac{1}{z}\big).\nn
\ee
Next we will show that the first integral vanishes as $\sqrt{z}\to \infty$. We begin with the case $\im (\sqrt{z})\to \infty$. 
For that purpose note that
\[
k(z)=\sqrt{z}+c+O(z^{-1/2}), \quad \textrm{as $z\to \infty$},
\]
for some constant $c\in \C$. Thus we compute
\begin{align*}
\Big| &\int_x^{\pm \infty}\frac{\psi_{q,\pm}(z,y)^2}{\psi_{q,\pm}(z,x)^2}\what{V}(y)dy\Big|
\leq C \int_x^{\pm \infty}\exp\big(\mp 2\im ( \sqrt{z})(y-x)\big)\big|\what{V}(y)\big|dy\\
&\leq C\int_{x}^{x+\varepsilon}\big|\what{V}(y)\big|dy+
C\cdot \exp\big(\mp2\im (\sqrt{z})\varepsilon\big)\int_{x+\varepsilon}^{\pm \infty}\big|\what{V}(y)\big|dy,
\end{align*}
such that the first integral can be made arbitrary small if $\varepsilon >0$ is small and the second integral vanishes as $\im (\sqrt{z})\to \infty$.

Otherwise, if $\re (\sqrt{z})\to \infty$, we use \eqref{decomppsipm} to rewrite the integral as
\[
\int_x^{\pm \infty } \left( \frac{\theta_{q,\pm}(z,y)^2}{\theta_{q,\pm}(z,x)^2} \what{V}(y) \exp\big(\mp 2\im ( \sqrt{z})(y-x)\big) \right)
\exp\big(\pm 2\I\re ( \sqrt{z})(y-x)\big) dy
\]
Since
\[
\left| \frac{\theta_{q,\pm}(z,y)^2}{\theta_{q,\pm}(z,x)^2} \what{V}(y) \exp\big(\mp 2\im ( \sqrt{z})(y-x)\big) \right| \le C |\what{V}(y)|
\]
the integral vanishes as $\re (\sqrt{z})\to \infty$ by a slight variation of the Riemann--Lebesgue lemma.

Hence we finally have
\be
\frac{\psi_\pm(z,x)}{\psi_{q,\pm}(z,x)}=1\pm \frac{\I }{2\sqrt{z}}\int_{x}^{\pm \infty}\what{V}(y)dy+o\big(\frac{1}{\sqrt{z}}\big)
\ee
as $z\to \infty$.
\end{proof}

\noindent
For later use we note the following immediate consequence

\begin{corollary} \label{corpisprime}
Under the assumptions of the previous theorem we have
\[
\lim_{x\to \pm \infty}\E^{\mp \I xk(z)}\Big( \dot{\psi}_\pm(z,x)\mp \I x \dot{k}(z)\psi_\pm(z,x)-
\dot{\psi}_{q,\pm}(z,x)\pm \I x \dot{k}(z)\psi_{q,\pm}(z,x)\Big)=0,
\]
where the dot denotes differentiation with respect to $z$.
\end{corollary}

\begin{proof}
Just differentiate \eqref{jost1} with respect to $z$, which is permissible by uniform
convergence on compact subsets of $\C\backslash \{E_j\}_{j=0}^{2g}$.
\end{proof}

\noindent
We remark that if we require our perturbation to satisfy the usual short-range
assumption as in \cite{BET}, \cite{F1,F2,F3} (i.e., the first moment is integrable, see \eqref{hypo2}), then we even have
$\E^{\mp \I xk(z)} (\dot{\psi}_\pm(z, x) - \dot{\psi}_{q,\pm}(z, x)) \to 0$.

From Theorem~\ref{thmjost} we obtain a complete characterization of the spectrum of $H$.

\begin{theorem} 
Assume \eqref{hypo}. Then $(H-z)^{-1}-(H_q-z)^{-1}$ is trace class. In particular,
we have $\sig_{ess}(H)=\sig(H_q)$ and the point spectrum of $H$ is confined to $\ol{\R\backslash\sig(H_q)}$.
Furthermore, the essential spectrum of $H$ is purely absolutely continuous except for possible
eigenvalues at the band edges.
\end{theorem}

\begin{proof}
That $(H-z)^{-1}-(H_q-z)^{-1}$ is trace class follows as in \cite[Lem.~9.34]{te} (cf.\ also \cite[Sect.~4]{krt}).
The fact that the essential spectrum is purely absolutely continuous follows from subordinacy theory
(\cite[Sect.~9.5]{te}) since the asymptotics of the Jost solutions imply that no solution is subordinate
inside the essential spectrum.
\end{proof}

\noindent
Note that \eqref{hypo} does neither exclude eigenvalues at the boundary of the essential spectrum nor an
infinite number of eigenvalues inside essential spectral gaps (see \cite{R-B}, \cite{krt3} or \cite{BET} for conditions
excluding these cases).

Our next result concerns the asymptotics of the Jost solutions at the {\em other side}.

\begin{lemma} \label{lemothers}
Assume \eqref{hypo}. Then the Jost solutions
$\psi_\pm(z, .)$, $z \in \C\backslash\sig(H)$, satisfy
\be                         \label{perturbed sol}
\lim_{x \rightarrow \mp \infty}
\big| \E^{\mp \I xk(z)} \big(\psi_\pm(z, x) - \alpha(z)\psi_{q,\pm}(z, x)\big)\big| = 0,
\ee
where
\be
\alpha(z) = \frac{W(\psi_-(z),\psi_+(z))}{W(\psi_{q,-}(z), \psi_{q,+}(z))} =
\frac {\prod_{j=1}^g(z - \mu_j)}{2\I \Rg{z}} W(\psi_-(z), \psi_+(z)).
\ee
\end{lemma}

\begin{proof}
Since $H$ and $H_q$ have the same form domain, the second resolvent equation (\cite[Lem.~6.30]{te})
for form perturbations implies
\[
G(z,x,x)- G_q(z,x,x) = \int_{-\infty}^\infty G(z,x,y) \what{V}(y) G_q(z,y,x) dy,
\]
where $ \what{V}= V-V_q$.
By \eqref{defGq} and 
\[
G(z,x,y)=\frac{1}{W( \psi_+,\psi_-)}\left\{ \begin{array}{cc}
\psi_+(z,x)\psi_-(z,y), & x\geq y, \\ \psi_+(z,y)\psi_-(z,x), & x\leq y,
\end{array}
\right.
\]
we obtain
\begin{align} \label{greendiff}
G(z,x,x)- G_q(z,x,x) =& \frac{\psi_{q,+}(z,x) \psi_+(z,x)}{W(z) W_q(z)}
\int_{-\infty}^x \what{V}(y) \psi_{q,-}(z,y) \psi_-(z,y)dy\nn \\
& + \frac{\psi_{q,-}(z,x) \psi_-(z,x)}{W(z) W_q(z)}
\int_x^\infty  \what{V}(y) \psi_{q,+}(z,y) \psi_+(z,y)dy,
\end{align}
where $W(z)=W( \psi_+,\psi_-)$ and $W_q(z)=W( \psi_{q,+},\psi_{q,-})$. Next, by \eqref{jost2}, note that
\be
\big| \psi_{q,\pm}(z,x)\big|\leq c_1 \E ^{\mp \varepsilon x},\quad 
\big| \psi_\pm(z,x)\big|\leq c_2 \E ^{\mp \varepsilon x},
\ee
as $x\to +\infty$, where $c_1$, $c_2$ denote some constants and $\varepsilon>0$ does only depend on $z$.  
Now one can show that the first term in \eqref{greendiff} tends to $0$ when $x\to+\infty$ using the same kind of argument as in the proof of Theorem~\ref{thmjost}. Similarly one then checks that the second term in \eqref{greendiff} tends to $0$ when $x\to-\infty$.
Thus
\[
\lim_{x\to\pm\infty} G(z,x,x)- G_q(z,x,x) =0
\]
and using
\[
G_q(z,x,x) = \frac{\psi_{q,-}(z,x) \psi_{q,+}(z,x)}{W(\psi_{q,-}(z),\psi_{q,+}(z))},
\qquad G(z,x,x) = \frac{\psi_-(z,x) \psi_+(z,x)}{W(\psi_-(z),\psi_+(z))}
\]
implies
\[
\lim_{x\to\pm\infty} \big(\psi_-(z,x)\psi_+(z,x) - \alpha(z) \psi_{q,-}(z,x)\psi_{q,+}(z,x)\big) =0,
\]
respectively,
\[
\lim_{x\to-\infty} \psi_{q,-}(z,x)\big(\psi_+(z,x) - \alpha(z) \psi_{q,+}(z,x) \big) =0,
\]
which is the claimed result.
\end{proof}

\noindent
To see the connection with scattering theory (see, e.g., \cite{BET}), we introduce the scattering relations
\begin{equation}\label{S2.16}
T(\lambda) \psi_\pm(\lambda,x) =\overline{\psi_\mp(\lambda,x)} +
R_\mp(\lambda)\psi_\mp(\lambda,x), \quad\lambda\in\sig(H_q),
\end{equation}
where the transmission and reflection coefficients are defined as
usual,
\begin{equation}
T(\lambda)= \frac{W(\overline{\psi_\pm(\lambda)}, \psi_\pm(\lambda))}{W(\psi_\mp(\lambda),
\psi_\pm(\lambda))}, \qquad
R_\pm(\lambda):= - \frac{W(\psi_\mp(\lambda),\overline{\psi_\pm(\lambda)})}
{W(\psi_\mp(\lambda), \psi_\pm(\lambda))},
\quad\lambda\in \sig(H_q).
\end{equation}
In particular, $\alpha(z)$ is just the inverse of the transmission coefficient $T(z)$.
It is holomorphic in $\C\backslash\sig(H_q)$ with simple zeros at the discrete eigenvalues of $H$.

\begin{corollary}
Assume \eqref{hypo}. Then we have
\be
T(z)=\exp \Big( - \int_{-\infty}^{+\infty}\big(m_\pm(z,x)-m_{q,\pm}(z,x)\big)dx\Big),
\ee
where 
\be \label{defweylm}
m_\pm(z,x)= \pm\frac{\psi_\pm'(z,x)}{\psi_\pm(z,x)}, \quad m_{q,\pm}(z,x)= \pm\frac{\psi_{q,\pm}'(z,x)}{\psi_{q,\pm}(z,x)}
\ee
are the Weyl--Titchmarsh functions.
Here the prime denotes differentiation with respect to $x$.
\end{corollary}

\begin{proof}
From the definition \eqref{defweylm} we get the following representations of the Jost and Baker--Akhiezer functions
\begin{align*}
\psi_\pm(z,x)=&\psi_\pm(z,x_0)\exp \big( \pm\int_{x_0}^xm_\pm(z,y)dy\big),\\
\psi_{q,\pm}(z,x)=&\psi_{q,\pm}(z,x_0)\exp \big( \pm\int_{x_0}^xm_{q,\pm}(z,y)dy\big),
\end{align*}
and thus
\begin{align*}
\frac{\psi_\pm(z,x)}{\psi_{q,\pm}(z,x)}=&\frac{\psi_\pm(z,x_0)}{\psi_{q,\pm}(z,x_0)}
\exp \big( \pm\int_{x_0}^x(m_\pm(z,y)-m_{q,\pm}(z,y))dy\big)\\
=& \exp \big( \pm\int_{\pm \infty}^x(m_\pm(z,y)-m_{q,\pm}(z,y))dy\big).
\end{align*}
Making use of that and \eqref{perturbed sol} we get
\[
\alpha(z)=\lim_{x \rightarrow \mp \infty}\frac{\psi_\pm(z,x)}{\psi_{q,\pm}(z,x)}=
\exp \big( \pm\int_{\pm \infty}^{\mp\infty}(m_\pm(z,y)-m_{q,\pm}(z,y))dy\big),
\]
which finishes the proof.
\end{proof}

\begin{corollary}\label{corasymal}
Assume \eqref{hypo}. Then $T(z)$ has the following asymptotic behavior
\be
T(z) = 1 + \frac{1}{2\I \sqrt{z}} \int_{-\infty}^{\infty} \big(V(y)-V_q(y)\big) dy + o(z^{-1/2}),
\ee
as $z\to \infty$.
\end{corollary}

\begin{proof}
Use \eqref{perturbed sol} and \eqref{B4jost}.
\end{proof}

\section{Connections with Krein's spectral shift theory and trace formulas}

To establish the connection with Krein's spectral shift theory we next
show:

\begin{lemma}
We have
\be
\frac{d}{dz} \alpha(z) = - \alpha(z) \int_{-\infty}^{+\infty} \big( G(z, x, x) - G_q(z, x, x)\big) dx,
\qquad z\in\C\backslash\sig(H),
\ee
where $G(z,x,y)$ and $G_q(z,x,y)$ are the Green's functions of $H$ and $H_q$,
respectively.
\end{lemma}

\begin{proof}
The Lagrange identity (\cite{te}, eq. (9.4)) implies
\be     \label{green 1}
W_x(\psi_+(z), \dot{\psi}_-(z)) - W_{y}(\psi_+(z), \dot{\psi}_-(z)) =  
\int_y^x \psi_+(z,r) \psi_-(z,r)dr, 
\ee 
hence the derivative of the Wronskian can be written as
\begin{align} \nn
&\frac{d}{dz}W(\psi_-(z), \psi_+(z)) = W_x(\dot{\psi}_-(z), \psi_+(z)) +
W_x(\psi_-(z), \dot{\psi}_+(z)) \\ \nn
& \qquad = W_y(\dot{\psi}_-(z), \psi_+(z)) + W_x(\psi_-(z), \dot{\psi}_+(z)) -
\int_y^x \psi_+(z,r)\psi_-(z,r)dr.
\end{align}
Using Corollary~\ref{corpisprime} and Lemma~\ref{lemothers} 
we have 
\bea \nn
W_y(\dot{\psi}_-(z), \psi_+(z)) &=& W_y(\dot{\psi}_- +\I \dot{k} y \psi_-, \psi_+) -\\
\nn
&& \I \dot{k} \big( y\, W(\psi_-, \psi_+) - \psi_-(z,y) \psi_+(z,y) \big)\\ \nn
&\to& \alpha W_{y} (\dot{\psi}_{q,-} + \I \dot{k} y \psi_{q,-}, \psi_{q,+}) -\\ \nn
&& \alpha \I \dot{k} \big( y\, W(\psi_{q,-}, \psi_{q,+}) -
\psi_{q,-}(z,y) \psi_{q,+}(z,y) \big)\\ \nn
&=& \alpha(z) W_y(\dot{\psi}_{q,-}(z), \psi_{q,+}(z))
\eea
as $y \rightarrow - \infty$. Similarly we obtain
\[
W_x(\psi_-(z), \dot{\psi}_+(z)) \to \alpha(z)
W_x(\psi_{q,-}(z), \dot{\psi}_{q,+}(z))
\]
as $x \rightarrow +\infty$
and again using \eqref{green 1} we have
\[
W_y(\dot{\psi}_{q,-}(z), \psi_{q,+}(z)) 
= W_x(\dot{\psi}_{q,-}(z), \psi_{q,+}(z))
+ \int_y^x \psi_{q,+}(z,r) \psi_{q,-}(z,r)dr.
\] 
Collecting terms we arrive at
\begin{align} \nn
\dot{W}(\psi_-(z), \psi_+(z)) = & 
- \int_{-\infty}^{+\infty} \Big( \psi_+(z, r) \psi_-(z, r) - 
\alpha(z) \psi_{q,+}(z, r) \psi_{q,-}(z, r) \Big)dr \\  \nn
& + \alpha(z) \dot{W}(\psi_{q,-}(z) \psi_{q,+}(z)).
\end{align}
Abbreviating $W_q=W(\psi_{q,-},\psi_{q,+})$ we now compute
\begin{align} \nn
\frac{d}{dz} \alpha(z) &= \frac{d}{dz} \Big( \frac{W}{W_q}\Big)
= - \frac{\dot{W}_q}{W_q^2} W + \frac{1}{W_q}
\Big( - \int_{-\infty}^{+\infty} \big( \psi_+  \psi_-  - \alpha  \psi_{q,+}  \psi_{q,-} \big)dr    
+ \alpha  \dot{W}_q\Big) \\ \nn
&= - \frac{1}{W_q} \int_{-\infty}^{+\infty} \Big( \psi_+ (z,r) \psi_-(z,r) - \alpha(z)  \psi_{q,+}(z,r)  \psi_{q,-}(z,r) \Big)dr, 
\end{align}
which finishes the proof.
\end{proof}

\noindent
Since $(H-z)^{-1} - (H_q-z)^{-1}$ is trace class with continuous integral kernel $G(z,x,x)-G_q(z,x,x)$,
we have (\cite{br})
\be
\tr\big((H-z)^{-1} - (H_q-z)^{-1} \big) = \int_{-\infty}^{+\infty} \big( G(z, x, x) - G_q(z, x, x)\big) dx,
\qquad z\in\C\backslash\sig(H),
\ee
and the last result can be rephrased as
\be\label{alrestr}
\frac{d}{dz} T(z) = T(z) \tr\big((H-z)^{-1} - (H_q-z)^{-1} \big),
\qquad z\in\C\backslash\sig(H),
\ee

As an immediate consequence we can establish the connection with Krein's spectral shift function (\cite{krein}).
We refer to \cite{yafams} for Krein's spectral shift theory in the case when only the resolvent difference is trace class;
which is the case needed here.

\begin{theorem}
The transmission coefficient $T(z)$ has the representation
\be\label{repalpha}
T(z) =\exp \Big(\int_\R \frac{\xi(\lambda)d\lambda}{\lambda - z} \Big),
\ee
where
\be
\xi(\lambda) = \frac{1}{\pi}\lim_{\epsilon \downarrow 0} 
\arg T(\lambda + \I\epsilon) 
\ee
is the spectral shift function of the pair $H$, $H_q$. Moreover,
$(V-V_q)^{1/2} (H_q-z)^{-1} |V-V_q|^{1/2}$ is trace class and $T(z)$
is the perturbation determinant of the pair $H$ and $H_q$:
\be
T(z) =\det \big(\id + (V-V_q)^{1/2} (H_q-z)^{-1} |V-V_q|^{1/2}\big).
\ee
If in addition $(V-V_q) (H_q-z)^{-1}$ is trace class we have
\be
T(z) =\det \big(\id + (V-V_q) (H_q-z)^{-1}\big).
\ee
\end{theorem}

\begin{proof}
The function $\im\log(T(z))$ is a bounded harmonic function in the upper half
plane and hence has a Poisson representation (cf.\ \cite{ko})
\[
\im\log(T(z)) =\int_{\R} \frac{y}{(x-\lambda)^2 + y^2} \xi(\lambda)d\lambda.
\quad z=x+\I y,
\]
Moreover, by $\xi(\lambda)=0$ for $\lam$ below the spectrum of $H$ and
$\xi(\lambda)= O(\lam^{-1/2})$ as $\lam\to+\infty$ (by Corollary~\ref{corasymal}) we obtain equality in \eqref{repalpha}
up to a real constant. The missing constants follows since both sides tend to $1$ as $z\to\infty$.
Moreover, combining \eqref{repalpha} with \eqref{alrestr} we see
\[
\tr\big((H-z)^{-1} - (H_q-z)^{-1} \big) = \int_{\R} 
\frac{\xi(\lambda)d\lambda}{(\lambda - z)^2},
\]
which shows that $\xi(\lambda)$ is the spectral shift function.

That $T(z)$ is the perturbation determinant is standard if $(V-V_q) (H_q-z)^{-1}$ is trace class
(see e.g.\ \cite{yafams}) for the slightly more general case when $(V-V_q)^{1/2} (H_q-z)^{-1} |V-V_q|^{1/2}$
is trace class we refer to \cite[Sect.~4]{GM}, \cite[Sect.~7]{glmz}. That this last condition holds will be shown in the
next lemma below.
\end{proof}

To following result needed in the previous proof is of independent interest.

\begin{lemma}
Assume \eqref{hypo}. Then $(V-V_q)^{1/2}(H_q-z)^{-1} |V-V_q|^{1/2}$ is trace class.
If we even have
\be\label{trbsc}
\|V-V_q\|_{2;1} = \sum_{n\in\Z} \left( \int_n^{n+1} |V(x)-V_q(x)|^2\right)^{1/2} < \infty,
\ee
then $(V-V_q) (H_q-z)^{-1}$ is trace class.
\end{lemma}

\begin{proof}
To see the first claim we begin with the fact \cite[Prop.~2.2]{sres} that $|V-V_q|^{1/2}(H_0-z)^{-1} |V-V_q|^{1/2}$ 
is trace class, where $H_0=-\frac{d^2}{dx^2}$. Let $z<0$ and set $A(z) = (V-V_q)^{1/2}(H_0-z)^{-1/2}$.
Then $A(z) A(z)^* = |V-V_q|^{1/2}(H_0-z)^{-1} |V-V_q|^{1/2}$ is trace class and thus $A(z)$ is Hilbert--Schmidt.
In fact, since $A(z) = A(z_0) (H_0-z_0)^{1/2}(H_0-z)^{-1/2}$ this holds for all $z\in\rho(H_0)$ and not just for $z<0$.
Hence, using $(H_q-z)^{-1} = (H_0-z)^{-1/2} C(z) (H_0-z)^{-1/2}$, where $C(z)$ is bounded
(cf.\ \cite[Thm.~6.25]{te}), we see $ |V-V_q|^{1/2}(H_q-z)^{-1} |V-V_q|^{1/2} =  A C(z) A^*$ which
establishes the claim.

The see the second claim we again begin with the fact \cite[Theorem~4.5]{str}  that \eqref{trbsc}
implies that $(V-V_q) (H_0-z)^{-1}$ is trace class. Now the second resolvent equation
$(H_q-z)^{-1} = (H_0-z)^{-1} - (H_0-z)^{-1}V_q (H_q-z)^{-1}$
establishes the claim since $V_q(H_q-z)^{-1}$ is bounded (cf.\ \cite[Sect.~9.7]{te}).
\end{proof}

Note that in the case $V_q=0$ \cite[Prop~4.7]{str} implies that the condition \eqref{trbsc} is optimal.
Moreover, the norm in \eqref{trbsc} dominates the $L^1$ norm, $\|V\|_1 \le \|V\|_{2;1}$ by the
Cauchy--Schwartz inequality, but the converse is of course not true (since \eqref{trbsc} forces the
function to be locally square integrable).

\section{The transmission coefficient}

Throughout this section we make the somewhat stronger assumption that
\be \label{hypo2}
\int_{-\infty}^{+\infty} (1+|x|)\big|V(x) - V_q(x)\big| dx <\infty
\ee
in order to ensure that there is only a finite number of eigenvalues in each gap \cite{R-B}.
Our aim is to reconstruct the transmission coefficient $T(z)$ from its boundary values
and its poles. To this end, recall that $T(z)$ is meromorphic in $\C\backslash\sig(H_q)$ with
simple poles at the eigenvalues $\rho_j$ of $H$. Moreover, for $z\in\sig(H_q)$ the
boundary values from the upper, respectively, lower, half plane exist and satisfy
$|T(z)|^2 = 1 - |R_\pm(z)|^2$, where $R_\pm(z)$ are the reflection coefficients defined in the
previous section.

In the case where $V_q=0$, this can be done via the classical Poisson--Jensen formula.
In the more general setting here, the reconstruction needs to be done on the
underlying Riemann surface. We essentially follow \cite{tag}, where the analog problem
for Jacobi operators was solved.

Denote by $\om_{p\, q}$ the normalized Abelian differential of the third kind with poles
at $p$ and $q$. Then the Blaschke factor is defined by
\be \label{Blaschke}
B(p,\rho)= \exp \Big( g(p,\rho) \Big) = \exp\Big(\int_{E_0}^p \om_{\rho\, \rho^*}\Big) =
\exp\Big(\int_{E(\rho)}^\rho \om_{p\, p^*}\Big), \quad \pi(\rho)\in\R,
\ee
where $E(\rho)$ is $E_0$ if $\rho<E_0$ and either $E_{2j-1}$ or $E_{2j}$ if
$\rho\in(E_{2j-1},E_{2j})$, $1\le j \le g$.
It is a multivalued function with a simple zero at $\rho$ and simple pole at $\rho^*$
satisfying $|B(p,\rho)|=1$, $p\in\partial\Pi_+$. It is real-valued for $\pi(p)\in(-\infty,E_0)$ and
satisfies
\be\label{eq:propblaschke}
B(E_0,\rho)=1 \quad\mbox{and}\quad
B(p^*,\rho) = B(p,\rho^*) = B(p,\rho)^{-1}.
\ee

Then we have

\begin{theorem}
The transmission coefficient is given by
\be \label{defd}
T(z,x) = \bigg( \prod_{j=1}^g B(p,\rho_j)^{-1}\bigg)\exp\Big( \frac{1}{2\pi\I}
\int_{\partial \Pi_+} \log (1-|R_\pm|^2) \om_{p\,p_{\infty}}\Big),
\quad p=(z,+),
\ee
where we set $R_\pm(p)= R_\pm(z)$ for $p=(z,+)$ and $R_\pm(p)=\ol{R_\pm(z)}$ for $p=(z,-)$.
\end{theorem}

\begin{proof}
Just literally follow the argument in \cite[Sect.~3]{tag}.
\end{proof}

\begin{remark}
A few remarks are in order:
\begin{enumerate}
\item
Using symmetry, $|R_\pm(p^*)|=|R_\pm(p)|$ for $p\in\partial\Pi_+$, of the integrand we can rewrite \eqref{defd} as
\be \label{defd2}
T(p,x) = \bigg( \prod_{j=1}^g \exp\Big(-\int_{E(\rho_j)}^{\rho_j} \om_{p\, p^*}\Big) \bigg)\exp\Big( \frac{1}{2\pi\I} \int_{\Sigma} \log (1-|R_\pm|^2) \om_{p\,p^*}\Big),
\ee
where the integral over $\Sigma$ is taken on the upper sheet.
\item
There exist explicit formulas for Abelian differentials of the third kind:
\be
\aligned
\om_{p q} &= \left( \frac{\Rgo + \Rg{p}}{2(\pi - \pi(p))} -
\frac{\Rgo + \Rg{q}}{2(\pi - \pi(q))} + P_{p q}(\pi) \right)
\frac{d\pi}{\Rgo},\nn \\
\om_{p p_{\infty}} &=\left( -\frac{\Rgo + \Rg{p}}{2(\pi - \pi(p))}+ P_{p p_{\infty}}(\pi) \right)
\frac{d\pi}{\Rgo},
\endaligned
\ee
where $P_{p q}(z)$, $P_{p p_{\infty}}(z)$ are polynomials of degree $g-1$ which have to be determined from
the normalization condition to have vanishing $a$-periods. In particular,
\[
\om_{p p^*} = \left( \frac{\Rg{p}}{\pi - \pi(p)} + P_{p p^*}(\pi) \right)
\frac{d\pi}{\Rgo}.
\]
\item
The function
\[
T(p) = \begin{cases} T(z), & p=(z,+),\\
T(z)^{-1}, & p=(z,-), \end{cases}
\]
solves the  following scalar meromorphic Riemann--Hilbert factorization problem:
\be \label{rhpptc}
\aligned
&T_+(p,x) = T_-(p,x) (1-|R(p)|^2), \quad p \in \partial \Pi_+,\\
&(T(p,x))= \dirhos -\dirho \\
&T(p_{\infty},x) = 1.
\endaligned
\ee
Here the subscripts in $T_\pm(p)$ denote the limits from $\Pi_\pm$, respectively.
Compare \cite{kt}, \cite{kt2}.
\end{enumerate}
\end{remark}

As was pointed out in \cite{tag}, this implies the following algebraic constraint on the
scattering data.

\begin{theorem}
The transmission coefficient $T$ defined via \eqref{defd} is single-valued if and
only if the eigenvalues $\rho_j$ and the reflection coefficients $R_\pm$
satisfy
\be \label{algcon}
\sum_j \int_{\rho_j^*}^{\rho_j} \zeta_\ell -
\frac{1}{2\pi\I} \int_{\partial \Pi_+}  \!\!\! \log(1 - |R_\pm|^2) \zeta_\ell \in \Z.
\ee
\end{theorem}

\section{Conserved quantities of the KdV hierarchy}

Finally we turn to solutions of the KdV hierarchy (see \cite{GH}). Let $V_q(x,t)$ be a finite-gap solution of
some equation in the KdV hierarchy, $\KdV_r(V_q(x,t))=0$, and let $V(x,t)$
be another solution, $\KdV_r(V(x,t))=0$, such that $V(.,t)-V_q(.,t)$ is Schwartz class for all $t\in\R$.
Existence of such solutions has been established only recently in \cite{EGT}.

Since the transmission coefficient $T(z,t)=T(z,0)\equiv T(z)$ is conserved (see
\cite{EGT} -- formally this follows from unitary invariance of the determinant),
conserved quantities of the KdV hierarchy can be obtained by computing
the asymptotic expansion at $\infty$.

To this end, we begin by recalling the following well-known asymptotics
for the Weyl $m$-functions in case of smooth potentials:

\begin{lemma} \label{lemdecomp}
Suppose $V(x) \in C^\infty(\R)$ is smooth.
The Weyl $m$-functions have the following asymptotic expansion for large $z$
\be\label{1.31}
m_\pm(z,x) \asymp \I \sqrt{z} \pm \sum_{n=1}^\infty \frac{\chi_n(x)}{(\pm 2\I\sqrt{z})^n},
\ee
with coefficients defined recursively via
\be\label{1.33}
\chi_1(x)=V(x),\quad\chi_{n+1}(x)=-\frac{\partial}
{\partial x}\chi_n(x) - \sum_{m=1}^{n-1}\chi_{n-m}(x)\chi_m(x).
\ee
\end{lemma}

The corresponding expansion coefficients associated with $V_q$ will be denoted by
$\chi_{q,m}(x)$. It is also known that the even coefficients are complete differentials  \cite{dj} and
the first few are explicitly given by
\begin{align*}
\chi_1(x) &= V(x),\\
\chi_2(x) &= -V'(x),\\
\chi_3(x) &= V''(x) - V(x)^2,\\
\chi_4(x) &= -V'''(x) + 4 V(x) V'(x),\\
\chi_5(x) &= V''''(x) - 6 V''(x) V(x) - 5 V'(x)^2 +2 V(x)^3.
\end{align*}

\begin{theorem}
Suppose $V(x)-V_q(x) \in \mathcal{S}(\R)$ is Schwartz. Then $\log T(z)$ has an asymptotic expansion
around $z=\infty$:
\[
\log T(z) \asymp \I \sqrt{z} \sum_{k=1}^\infty \frac{\tau_k}{z^k}.
\]
The quantities $\tau_k$ are given by
\be\label{eqtauk}
\tau_k=\int_{-\infty}^{\infty}\frac{\chi_{2k-1}(x) - \chi_{q,2k-1}(x)}{(-1)^k2^{2k-1}}dx
\ee
and are conserved quantities for the KdV hierarchy. Explicitly,
\begin{align*} \nn
\tau_1 =& -\frac{1}{2}\int_{-\infty}^{\infty}\big( V(x)-V_q(x)\big)dx,\\
\tau_2 =& -\frac{1}{8}\int_{-\infty}^{\infty}\big(V^2(x) - V_q^2(x)\big)dx,\\
\tau_3 =& -\frac{1}{32}\int_{-\infty}^{\infty}\big(2 V^3(x) - 5 V_x^2(x) - 6 V_{xx}(x) V(x)\\
& \qquad -2 V_q^3(x) + 5 V_{q,x}^2(x) + 6 V_{q,xx}(x) V_q(x)\big)dx,
\end{align*}
etc.
\end{theorem}

\begin{proof}
Represent the Jost solutions in the form
\be\label{Jost1}
\psi_\pm(z,x)=\psi_{q,\pm}(z,x)\exp\left(\mp\int_x^{\pm\infty}
\big( m_\pm(z,y) - m_{q,\pm}(z,y)\big) dy \right).
\ee
Then iterating the Volterra integral equations \eqref{Volterrafrac} one sees that $\psi_\pm(z,x)$
have an asymptotic expansion uniformly with respect to $x$ and given by
\be\label{kappa}
\log \frac{\psi_\pm(z,x)}{\psi_{q,\pm}(z,x)} \asymp - \sum_{n=1}^\infty\frac{1}{(\pm 2\I\sqrt{z})^n} \int_x^{\pm\infty}
\big(\chi_n(y) - \chi_{q,n}(y)\big) dy.
\ee
Then, letting $x\to\mp\infty$ using \eqref{perturbed sol} yields \eqref{eqtauk}. In particular, equality of the plus and
minus cases shows that all even expansion coefficients must vanish (which alternatively also follows from the fact that
the even expansion coefficients are complete differentials).
\end{proof}

\begin{theorem}
Consider the expansion coefficients $\tau_k$ of $\log T(z)$ defined in \eqref{eqtauk}.
Then the following trace formulas are valid:
\be
\tau_k=2\I \sum_{j=1}^g\int_{E(\rho_j)}^{\rho_j}\om_{p_{\infty},2k-2}-\frac{1}{\pi }\int_{\Sigma}\log |T|^2\om_{p_{\infty},2k-2},
\ee
where $\om_{p_\infty,k}$ is the Abelian differential of the second kind with a pole of order $k+2$ at $p_{\infty}$.
\end{theorem}

\begin{proof}
From $\frac{d^k}{dz^k}\om_{p\, E_0}=k!\om_{p_{\infty},k-1}$ we get that
\begin{align*}
\om_{p\, E_0}&=\om_{p_{\infty}\, E_0}+\sum_{k=1}^{\infty}\zeta^k \om_{p_{\infty},k-1}, \qquad \zeta=z^{-1/2},\\
\om_{p^*\, E_0}&=\om_{p_{\infty}\, E_0}+\sum_{k=1}^{\infty}\zeta^k \om_{p_{\infty},k-1}, \qquad \zeta=-z^{-1/2}.
\end{align*}
Using this it follows
\[
\om_{p\, p^*}=\om_{p\, E_0}-\om_{p^*\, E_0}=2\sum_{k=1}^{\infty}\om_{p_{\infty},2k-2}\zeta^{2k-1},\quad\zeta=z^{-1/2}.
\]
Hence we have
\begin{align*}
&-\sum_{j=1}^g\int_{E(\rho_j)}^{\rho_j}\om_{p\, p^*}+\frac{1}{2\pi \I}\int_{\Sigma}\log |T|^2\om_{p\, p^*}=\\
&-\sum_{j=1}^g\int_{E(\rho_j)}^{\rho_j}2\sum_{k=1}^{\infty}\zeta^{2k-1} \om_{p_{\infty},2k-2}+\frac{1}{\pi \I}\int_{\Sigma}\log |T|^2\sum_{k=1}^{\infty}\zeta^{2k-1} \om_{p_{\infty},2k-2}.
\end{align*}
Thus, since $|T|^2 = 1 - |R_\pm|^2$ and $R_\pm(\lambda)$ decays faster than any polynomial as $\lambda\to\infty$
\cite{EGT}, one obtains
\[
\log T(z) \asymp -\sum_{k=1}^{\infty}\big( 2\sum_{j=1}^{\infty}\int_{E(\rho_j)}^{\rho_j} \om_{p_{\infty},2k-2}-\frac{1}{\pi \I}\int_{\Sigma}\log |T|^2\om_{p_{\infty},2k-2}\big)\zeta^{2k-1},
\]
where $\zeta=z^{-1/2}$ denotes the local coordinate at $z=\infty$.
\end{proof}

\begin{remark}
The differentials $\om_{p_{\infty},2k-2}$, $k=1,2,\dots$, are explicitly given by
\be
\om_{p_\infty,2k-2}=\left( \frac{\pi^{g+k-1}}{\Rgo}+P_k(\pi)\right)d\pi.
\ee
Here $P_k(\pi)$ is a polynomial of degree $g+k-2$ which has to be chosen such that all $a$-periods vanish.
\end{remark}

\section*{Acknowledgments}

G.T. would like to thank all organizers of the international conference on Operator Theory, Analysis and Mathematical Physics (OTAMP), Bedlewo, June 2008, for their kind invitation and the stimulating atmosphere during the meeting.
We are indebted to Iryna Egorova and Fritz Gesztesy for several helpful discussions.

\end{document}